\documentclass[11pt, letterpaper]{article}
\usepackage[utf8]{inputenc}
\usepackage{amssymb, amsfonts, amsmath, amsthm, mathtools}
\usepackage[margin=1.2in]{geometry}
\usepackage{enumitem}
\usepackage{ifthen}
\usepackage[none]{hyphenat}
\usepackage{theoremref}
\usepackage[usenames,dvipsnames]{xcolor}
\usepackage{tikz}
\usepackage{color}

\newtheoremstyle{mystyle}
{1.2em} 
{0.5em} 
    {\itshape} 
    {} 
    {\bfseries} 
    {.} 
    {.5em} 
    {} 
\theoremstyle{mystyle}
\newtheorem{thm}{Theorem}[section]
\theoremstyle{mystyle}
\newtheorem{prop}[thm]{Proposition}
\theoremstyle{mystyle}
\newtheorem{lemma}[thm]{Lemma}
\theoremstyle{mystyle}
\newtheorem{obs}[thm]{Observation}
\theoremstyle{mystyle}

\theoremstyle{mystyle}
\newtheorem{cor}[thm]{Corollary}
\theoremstyle{mystyle}

\theoremstyle{mystyle}
\newtheorem{que}{Question}
\theoremstyle{mystyle}
\newtheorem{conj}{Conjecture}

\input{tcilatex}

\definecolor{lgray}{gray}{0.95}
\definecolor{mgray}{gray}{0.40}
\usetikzlibrary{fit,positioning,calc}
\tikzstyle{std}=[ circle, draw=black,fill=black, inner sep=0pt, minimum size=2mm]
\tikzstyle{bred}=[circle, draw=black,fill=red,thick,  inner sep=2pt, minimum size=2mm]
\tikzstyle{bgreen}=[ circle, draw=black,fill=green,thick,  inner sep=2pt, minimum size=2.5mm]
\tikzstyle{sqRed}=[rectangle, draw=black,fill=red,thick,  inner sep=2pt, minimum size=2.5mm]
\tikzstyle{trEdge}=[color=black, line width = 1pt, style=dotted]

\begin{document}

\title{\textbf{Total Roman Domination Edge-Supercritical and Edge-Removal-Supercritical Graphs}}
\author{C. M. Mynhardt\thanks{%
Supported by a Discovery Grant from the Natural Sciences and Engineering
Research Council of Canada.}, S. E. A. Ogden\thanks{%
Supported by a Jamie Cassels Undergraduate Research Award from the University of
Victoria.} \\
%EndAName
Department of Mathematics and Statistics\\
University of Victoria\\
Victoria, BC, \textsc{Canada}\\
{\small kieka@uvic.ca, sogden@uvic.ca}}
\date{}
\maketitle

\begin{abstract}
A total Roman dominating function on a graph $G$ is a function $f:V(G)\rightarrow \{0,1,2\}$ such that every vertex $v$ with $f(v)=0$ is adjacent to some vertex $u$ with $f(u)=2$, and the subgraph of $G$ induced by the set of all vertices $w$ such that $f(w)>0$ has no isolated vertices. The weight of $f$ is $\Sigma _{v\in V(G)}f(v)$. The total Roman domination number $\gamma _{tR}(G)$ is the minimum weight of a total Roman dominating function on $G$. A graph $G$ is $k$-$\gamma _{tR}$-edge-critical if $\gamma_{tR}(G+e)<\gamma _{tR}(G)=k$ for every edge $e\in E(\overline{G})\neq\emptyset $, and $k$-$\gamma _{tR}$-edge-supercritical if it is $k$-$\gamma_{tR}$-edge-critical and $\gamma _{tR}(G+e)=\gamma _{tR}(G)-2$ for every edge $e\in E(\overline{G})\neq \emptyset $. A graph $G$ is $k$-$\gamma _{tR}$-edge-stable if $\gamma_{tR}(G+e)=\gamma _{tR}(G)=k$ for every edge $e\in E(\overline{G})$ or $E(\overline{G})=\emptyset$. For an edge $e\in E(G)$ incident with a degree $1$ vertex, we define $\gamma_{tR}(G-e)=\infty$. A graph $G$ is $k$-$\gamma _{tR}$-edge-removal-critical if $\gamma_{tR}(G-e)>\gamma _{tR}(G)=k$ for every edge $e\in E(G)$, and $k$-$\gamma _{tR}$-edge-removal-supercritical if it is $k$-$\gamma_{tR}$-edge-removal-critical and $\gamma _{tR}(G-e)\geq\gamma _{tR}(G)+2$ for every edge $e\in E(G)$. A graph $G$ is $k$-$\gamma_{tR}$-edge-removal-stable if $\gamma_{tR}(G-e)=\gamma_{tR}(G)=k$ for every edge $e\in E(G)$. We investigate connected $\gamma_{tR}$-edge-supercritical graphs and exhibit infinite classes of such graphs. In addition, we characterize $\gamma_{tR}$-edge-removal-critical and $\gamma_{tR}$-edge-removal-supercritical graphs. Furthermore, we present a connection between $k$-$\gamma_{tR}$-edge-removal-supercritical and $k$-$\gamma_{tR}$-edge-stable graphs, and similarly between $k$-$\gamma_{tR}$-edge-supercritical and $k$-$\gamma_{tR}$-edge-removal-stable graphs.\end{abstract}

\noindent\textbf{Keywords:} total domination; total Roman domination; total Roman domination edge-critical graphs; total Roman domination edge-supercritical graphs; total Roman domination edge-stable graphs; total Roman domination edge-removal-critical graphs; total Roman domination edge-removal-supercritical graphs; total Roman domination edge-removal-stable graphs

\noindent\textbf{AMS Subject Classification Number 2010:} 05C69

\section{Introduction}
\label{Sec: Intro}

We consider the behaviour of the total Roman domination number of a graph $G$ upon the addition or removal of edges to and from $G$. A \emph{dominating set} $S$ in a graph $G$ is a set of vertices such that every vertex in $V(G)-S$ is adjacent to at least one vertex in $S$. The \emph{domination number} $\gamma (G)$ is the cardinality of a minimum dominating set in $G$. A \emph{total dominating set} $S$ (abbreviated by \emph{TD-set}) in a graph $G$ with no isolated vertices is a set of vertices such that every vertex in $V(G)$ is adjacent to at least one vertex in $S$. The \emph{total domination number} $\gamma _{t}(G)$ (abbreviated by \emph{TD-number}) is the cardinality of a minimum total dominating set in $G$. For $S\subseteq V(G)$ and a function $f:S\rightarrow \mathbb{R}$, define $f(S)=\Sigma _{s\in S}f(s)$. A \emph{Roman dominating function} (abbreviated by \emph{RD-function}) on a graph $G$ is a function $f:V(G)\rightarrow\{0,1,2\}$ such that every vertex $v$ with $f(v)=0$ is adjacent to some vertex $u$ with $f(u)=2$. The \emph{weight} of $f$, denoted by $\omega (f)$, is defined as $f(V(G))$. The \emph{Roman domination number} $\gamma _{R}(G)$ (abbreviated by \emph{RD-number}) is defined as $\min \{\omega (f):f\text{ is an RD-function on }G\}$. For an RD-function $f$, let $V_{f}^{i}=\{v\in V(G):f(v)=i\}$ and $V_{f}^{+}=V_{f}^{1}\cup V_{f}^{2}$. Thus, we can uniquely express an RD-function $f$ as $f=(V_f^0,V_f^1,V_f^2)$. 

As defined by Chang and Liu \cite{ChangL}, a \emph{total Roman dominating function} (abbreviated by \emph{TRD-function}) on a graph $G$ with no isolated vertices is a Roman dominating function with the additional condition that $G[V_{f}^{+}]$ has no isolated vertices. The \emph{total Roman domination number} $\gamma _{tR}(G)$ (abbreviated by \emph{TRD-number}) is the minimum weight of a TRD-function on $G$; that is, $\gamma _{tR}(G)=\min \{\omega (f):f\text{ is a TRD-function on }G\}$. A TRD-function $f$ such that $\omega (f)=\gamma _{tR}(G)$ is called a $\gamma_{tR}(G)$-\emph{function}, or a $\gamma _{tR}$-\emph{function} if the graph $G$ is clear from the context; $\gamma_{R}$\emph{-functions} are defined analogously. Total Roman domination was also studied by Ahangar, Henning, Samodivkin and Yero \cite{AHSY}. 

The addition of an edge to a graph has the potential to change its total domination or total Roman domination number. Van der Merwe, Mynhardt and Haynes \cite{MMH} studied $\gamma_{t}$-\emph{edge-critical graphs}, that is, graphs $G$ for which $\gamma _{t}(G+e)<\gamma_{t}(G)$ for each $e\in E(\overline{G})$ and $E(\overline{G})\neq\emptyset $. Similarly, Lampman, Mynhardt and Ogden \cite{LMO} defined an edge $e\in E(\overline{G})$ to be \emph{critical} with respect to total Roman domination (abbreviated \emph{TRD-critical}) if $\gamma _{tR}(G+e)<\gamma _{tR}(G)$. An edge $e\in E(\overline{G})$ is \emph{supercritical} with respect to total Roman domination (abbreviated \emph{TRD-supercritical}) if $\gamma _{tR}(G+e)\leq \gamma_{tR}(G)-2$. A graph $G$ with no isolated vertices is \emph{total Roman domination edge-critical}, or simply $\gamma_{tR}$\emph{-edge-critical}, if every edge $e\in E(\overline{G})\neq\emptyset$ is TRD-critical. We say that $G$ is $k$\emph{-}$\gamma _{tR}$\emph{-edge-critical} if $\gamma _{tR}(G)=k$ and $G$ is $\gamma _{tR}$-edge-critical. Similarly, if every edge $e\in E(\overline{G})\neq\emptyset$ is TRD-supercritical, then $G$ is $\gamma _{tR}$\emph{-edge-supercritical}; $\gamma _{t}$\emph{-edge-supercritical} graphs are defined analogously. An edge $e\in E(\overline{G})$ is \emph{stable} with respect to total Roman domination (abbreviated \emph{TRD-stable}) if $\gamma_{tR}(G+e)=\gamma _{tR}(G)$. If every edge $e\in E(\overline{G})$ is TRD-stable, or if $E(\overline{G})=\emptyset$, we say that $G$ is $\gamma_{tR}$\emph{-edge-stable}.

The removal of an edge from a graph $G$ also has the potential to change its total domination or total Roman domination number. Desormeaux, Haynes and Henning \cite{DHH} studied $\gamma_{t}$-\emph{edge-removal-critical graphs}, that is, graphs $G$ for which $\gamma _{t}(G-e)>\gamma_{t}(G)$ for each $e\in E(G)$. We consider the same concept for total Roman domination. An edge $e\in E(G)$ is \emph{removal-critical} with respect to total Roman domination (abbreviated \emph{TRD-ER-critical}) if $\gamma_{tR}(G)<\gamma_{tR}(G-e)$. We say that an edge $e\in E(G)$ is \emph{removal-supercritical} with respect to total Roman domination (abbreviated \emph{TRD-ER-supercritical}) if $\gamma_{tR}(G)+2\leq\gamma _{tR}(G-e)$. Note that the removal of an edge $e\in E(G)$ incident with a degree $1$ vertex would result in $G-e$ containing an isolated vertex. For such an edge $e\in E(G)$, Desormeaux et al. \cite{DHH} defined $\gamma_{t}(G-e)=\infty$. Likewise, we define $\gamma_{tR}(G-e)=\infty$ when $e\in E(G)$ is an edge incident with a degree $1$ vertex. Furthermore, we define $E_P(G)\subseteq E(G)$ to be the set of edges in $G$ which are not incident with a degree $1$ vertex; that is, the set of edges $e$ such that $\gamma_{tR}(G-e)<\infty$. Hence every edge $e\in E(G)-E_P(G)$ is TRD-ER-supercritical. A graph $G$ with no isolated vertices is \emph{total Roman domination edge-removal-critical}, or simply $\gamma_{tR}$\emph{-ER-critical}, if every edge $e\in E(G)$ is TRD-ER-critical. We say that $G$ is $k$\emph{-}$\gamma _{tR}$\emph{-ER-critical} if $\gamma _{tR}(G)=k$ and $G$ is $\gamma _{tR}$-ER-critical. Similarly, if every edge $e\in E(G)$ is TRD-ER-supercritical, then $G$ is $\gamma _{tR}$\emph{-ER-supercritical}; $\gamma _{t}$\emph{-ER-supercritical} graphs are defined analogously. An edge $e\in E(G)$ is \emph{removal-stable} with respect to total Roman domination (abbreviated \emph{TRD-ER-stable}) if $\gamma_{tR}(G)=\gamma_{tR}(G-e)$. If every edge $e\in E(G)$ is TRD-ER-stable, we say that $G$ is $\gamma_{tR}$\emph{-edge-removal-stable}, or simply $\gamma_{tR}$\emph{-ER-stable}.

We refer the reader to the well-known books \cite{CL} and \cite{HHS} for graph theory concepts not defined here. Frequently used or lesser known concepts are defined where needed.

We begin with some previous results on the total domination and total Roman domination numbers of a graph in Section \ref{Sec: Prelim}, and $\gamma_{tR}$-edge-critical graphs in Section \ref{Sec: TRD crit}. In Section \ref{Sec: Super}, we investigate the existence of connected $\gamma_{tR}$-edge-supercritical graphs and demonstrate that each such graph contains a cycle. After characterizing $5$-$\gamma _{tR}$-edge-critical graphs in Section \ref{Sec: tR=5}, we investigate $6$-$\gamma_{tR}$-edge-supercritical graphs in Section \ref{Sec: 6-super}. In Section \ref{Sec: ER-critical}, we characterize $\gamma_{tR}$-ER-critical graphs. A similar characterization of $\gamma_{tR}$-ER-supercritical graphs is presented in Section \ref{Sec: ER-supercritical}, where we also note that every $\gamma_{tR}$-ER-supercritical graph is $\gamma_{tR}$-edge-stable. The analogous result for $\gamma_{tR}$-edge-supercritical and $\gamma_{tR}$-ER-stable graphs is given in Section \ref{Sec: ER-stable}. We conclude in Section \ref{Sec: Future} with ideas for future
research.

\section{Preliminaries}
\label{Sec: Prelim}

Before investigating $\gamma _{tR}$-edge-critical and $\gamma _{tR}$-ER-critical graphs, we present some basic results relating the domination, total domination, and total Roman domination numbers of a graph. Our first result is a direct corollary to Observation 6.42 and Theorem 6.47 in \cite{HHS}, and provides bounds on the total domination number of a graph $G$ in terms of its domination number.

\begin{prop}
\thlabel{TD vs D} \emph{\cite{HHS}} For a graph $G$ with no isolated vertices, $\gamma(G)\leq\gamma_{t}(G)\leq2\gamma(G)$. 
\end{prop}

\vspace{0mm}

As noted in Section \ref{Sec: Intro}, total Roman domination was studied by Ahangar et al. \cite{AHSY}. There, they provided two results which bound the total Roman domination number of a graph in terms of its domination number and total domination number, respectively. Note the similarities between the bounds in \thref{TD vs D,TRD vs TD}.

\begin{prop}
\thlabel{TRD vs D} \emph{\cite{AHSY}} For a graph $G$ with no isolated vertices, $2\gamma(G)\leq\gamma_{tR}(G)\leq3\gamma(G)$. 
\end{prop}

\begin{prop}
\thlabel{TRD vs TD} \emph{\cite{AHSY}} If $G$ is a graph with no isolated vertices, then $\gamma _{t}(G)\leq \gamma _{tR}(G)\leq 2\gamma _{t}(G)$. Furthermore, $\gamma _{tR}(G)=\gamma _{t}(G)$ if and only if $G$ is the disjoint union of copies of $K_2$. 
\end{prop}

\vspace{0mm}

Note that \thref{TRD vs TD} characterizes the graphs $G$ for which $\gamma_{tR}(G)=\gamma_{t}(G)$. Ahangar et al. \cite{AHSY} also characterized the graphs which nearly attain the lower bound in \thref{TRD vs TD}; that is, the graphs $G$ for which $\gamma_{tR}(G)=\gamma_{t}(G)+1$.

\begin{prop}
\thlabel{tR=t+1} \emph{\cite{AHSY}} Let $G$ be a connected graph of order $n\geq3$. Then $\gamma _{tR}(G)=\gamma _{t}(G)+1$ if and only if $\Delta(G)=n-1$, that is, $G$ has a universal vertex.
\end{prop}

\vspace{0mm}

We now consider the graphs with the smallest possible TRD-number, namely $3$, which were characterized by Lampman et al. \cite{LMO}.

\begin{prop}
\thlabel{tR=3} \emph{\cite{LMO}} For a graph $G$ of order $n\geq 3$ with no isolated vertices, $\gamma _{tR}(G)=3$ if and only if $\Delta (G)=n-1$, that is, $G$ has a universal vertex.
\end{prop}

\vspace{0mm}

When combined with \thref{tR=t+1}, \thref{tR=3} implies that, for a connected graph $G$ of order $n\geq 3$, $\gamma _{tR}(G)=\gamma _{t}(G)+1$ if and only if $\gamma _{tR}(G)=3$. This result provides a tighter lower bound on the TRD-number of a connected graph with no universal vertex with respect to its TD-number.

\begin{obs}
\thlabel{tR>t+1} If $G$ is a connected graph of order $n\geq 3$ such that $\Delta (G)\leq n-2$, then $\gamma _{t}(G)+2\leq \gamma _{tR}(G)\leq 2\gamma_{t}(G)$.
\end{obs}

\vspace{0mm}

Lampman et al. \cite{LMO} also provided an alternate characterization of the graphs $G$ with total Roman domination number $3$, as well as a characterization of the graphs $G$ with total Roman domination number $4$, in terms of the domination and total domination numbers of the graph.

\begin{prop}
\thlabel{tR=34} \emph{\cite{LMO}} If $G$ is a connected graph of order $n\geq 3$, then $\gamma_{tR}(G)\in \{3,4\}$ if and only if $\gamma_{t}(G)=2$. Moreover, $\gamma(G)=1$ when $\gamma_{tR}(G)=3$, and $\gamma(G)=2$ when $\gamma_{tR}(G)=4$.
\end{prop}

\section{$\gamma_{tR}$-Edge-critical graphs}
\label{Sec: TRD crit}

As noted in Section \ref{Sec: Intro}, the addition of an edge to a graph has the potential to change its total domination or total Roman domination number. Van der Merwe et al. \cite{MMH} studied this effect with respect to the total domination number, providing bounds on the total domination number of the graph $G+e$, where $e\in E(\overline{G})$, in terms of the total domination number of $G$. 

\begin{prop}
\thlabel{t bounds} \emph{\cite{MMH}} For a graph $G$ with no isolated vertices, if $uv\in E(\overline{G})$, then $\gamma_{t}(G)-2\leq \gamma _{t}(G+uv)\leq \gamma _{t}(G)$.
\end{prop}

\vspace{0mm}

These bounds also hold with respect to the total Roman domination number of the graph $G+e$ obtained by adding an edge $e\in E(\overline{G})$ to $G$, as shown by Lampman et al. \cite{LMO}.

\begin{prop}
\thlabel{tR bounds}  \emph{\cite{LMO}}
Given a graph $G$ with no isolated vertices, if ${uv\in E(\overline{G})}$, then $\gamma_{tR}(G)-2\leq \gamma _{tR}(G+uv)\leq \gamma _{tR}(G)$.
\end{prop}

\vspace{0mm}

For any edge $uv\in E(G)$, there are $3^2=9$ ways for a TRD-function $f$ to assign the values in $\{0,1,2\}$ to $u$ and $v$. However, the following observation restricts the possible values assigned to a degree $1$ vertex and its unique neighbour when $f$ is a $\gamma_{tR}(G)$-function. Note that, for a graph $G$ and a vertex $v\in V(G)$, the \emph{open neighbourhood} of $v$ in $G$ is $N_{G}(v)=\{u\in V(G):uv\in E(G)\}$, and the \emph{closed neighbourhood} of $v$ in $G$ is $N_{G}[v]=N_{G}(v)\cup \{v\}$.

\begin{obs}
\thlabel{pendant edge} For a graph $G$ with no isolated vertices, if $\text{deg}(u)=1$ and $N_G(u)=\{v\}$, then, for any $\gamma_{tR}(G)$-function $f$, either $f(u)=f(v)=1$, or $f(v)=2$ and $f(u)\in\{0,1\}$. Furthermore, there exists a $\gamma_{tR}(G)$-function $g$ such that $\{g(u),g(v)\}\neq\{1,2\}$. 
\end{obs}

\vspace{0mm}

Similarly, Lampman et al. \cite{LMO} provided a result restricting the possible values assigned to the vertices of a TRD-critical edge $uv$ by a $\gamma _{tR}$-function $f$ on $G+uv$. We mildly abuse set-theoretic notation by denoting the case where $f(u)=f(v)=i$ for $i\in\{0,1,2\}$ by $\{f(u),f(v)\}=\{i,i\}$.

\begin{prop}
\thlabel{set added edge} 
\emph{\cite{LMO}} Given a graph $G$ with no isolated vertices, if $uv\in E(\overline{G})$ is a TRD-critical edge and $f$ is a $\gamma _{tR}(G+uv)$-function, then $\{f(u),f(v)\}\in\{\{2,2\},\{2,1\},\{2,0\},\{1,1\}\}$. If, in addition, $\text{deg}(u)=\text{deg}(v)=1$, then there exists a $\gamma_{tR}(G+uv)$-function $f$ such that $f(u)=f(v)=1$.
\end{prop} 

\vspace{0mm}

We now consider $\gamma_{tR}$-edge-critical graphs. Recall that a graph $G$ with no isolated vertices is $\gamma _{tR}$-edge-critical if $\gamma_{tR}(G+e)<\gamma_{tR}(G)$ for every edge $e\in E(\overline{G})\neq\emptyset$. For a graph $G\neq K_{2}$, the unique neighbour of an end-vertex of $G$ is called its \emph{support vertex}. In this case, the end-vertex is referred to as a \emph{pendant vertex}, and the edge incident with it a \emph{pendant edge}. An \emph{endpath} in a graph $G$ is a path from a vertex $v$, where $\text{deg}(v)\geq3$, to a pendant vertex, such that all of the internal vertices of the path have degree $2$. We begin with some results from \cite{LMO} which provide necessary conditions for a graph $G$ to be $\gamma_{tR}$-edge-critical.

\begin{prop}
\thlabel{end deg 3} 
\emph{\cite{LMO}} For a graph $G$ with no isolated vertices, if $G$ has a pendant vertex $w$ with support vertex $x$ such that $G[N(x)-\{w\}]$ is not complete, then $G$ is not $\gamma _{tR}$-edge-critical. 
\end{prop}

\begin{prop}
\thlabel{no long legs} 
\emph{\cite{LMO}} For a graph $G$ with no isolated vertices, if $G$ has two endpaths $%
v_{0}, v_{1}, ..., v_{k}$ and $u_{0}, u_{1}, ..., u_{m}$, where $k,m\geq 3$ and $v_{k}$
and $u_{m}$ are pendant vertices, then $G$ is not $\gamma _{tR}$-edge-critical.
\end{prop}

\vspace{0mm}

We conclude this section by considering the graphs $G$ which have the largest TRD-number, namely $|V(G)|$. A \emph{subdivided star} is a tree obtained from a star on at least three vertices by subdividing each edge exactly once. A \emph{double star} is a tree obtained from two disjoint non-trivial stars by joining the two central vertices (choosing either central vertex in the case of $K_{2}$). The \emph{corona} $\func{cor}(G)$ (sometimes denoted by $G\circ K_{1}$) of $G$ is obtained by joining each vertex of $G$ to a new end-vertex.

Connected graphs $G$ for which $\gamma _{tR}(G)=|V(G)|$ were characterized in \cite{AHSY}. There, Ahangar et al. defined $\mathcal{G}$ as the family of connected graphs obtained from a $4$-cycle $v_{1},v_{2},v_{3},v_{4},v_{1}$ by adding $k_{1}+k_{2}\geq 1$ vertex-disjoint paths $P_{2}$, and joining $v_{i}$ to an end-vertex of $k_{i}$ such paths, for $i\in \{1,2\}$. Note that possibly $k_{1}=0$ or $k_{2}=0$. Furthermore, they defined $\mathcal{H}$ to be the family of graphs obtained from a double star by subdividing each pendant edge once and the non-pendant edge $r\geq 0$ times.

\begin{prop}
\thlabel{tR=n} \emph{\cite{AHSY}} If $G$ is a connected graph of order $%
n\geq 2$, then $\gamma _{tR}(G)=n$ if and only if one of the following
holds. 
\vspace{-2mm}
\setlist{nolistsep}
\begin{enumerate}
\item[$(i)$] $G$ is a path or a cycle;

\item[$(ii)$] $G$ is the corona of a graph;

\item[$(iii)$] $G$ is a subdivided star;

\item[$(iv)$] $G\in \mathcal{G}\cup \mathcal{H}$.
\end{enumerate}
\end{prop}

\vspace{0mm}

Lampman et al. \cite{LMO} used this result to characterize the connected graphs of order $n\geq4$ which are $n$-$\gamma _{tR}$-edge-critical. For $r\geq 0$, they defined $\mathcal{H}_{r}\subseteq \mathcal{H}$ as the family of graphs in $\mathcal{H}$ where the non-pendant edge was subdivided $r$ times. 

\begin{prop}
\thlabel{n edge-crit} \emph{\cite{LMO}} A connected graph $G$ of order $n\geq 4$ is $n$-$\gamma _{tR}$%
-edge-critical if and only if $G$ is one of the following graphs: %
\vspace{-2mm}
\setlist{nolistsep}
\begin{enumerate}
\item[$(i)$] $C_{n}$, $n\geq 4$;

\item[$(ii)$] $\func{cor}(K_{r})$, $r\geq 3$;

\item[$(iii)$] a subdivided star of order $n\geq 7$;

\item[$(iv)$] $G\in \mathcal{G}$;

\item[$(v)$] $G\in \mathcal{H}-\mathcal{H}_{0}-\mathcal{H}_{2}$.
\end{enumerate}
\end{prop}

\section{$\gamma_{tR}$-Edge-supercritical graphs}
\label{Sec: Super}

We now consider $\gamma _{t}$-edge-supercritical and $\gamma_{tR}$-edge-supercritical graphs. Note that, by \thref{t bounds}, a graph $G$ with no isolated vertices is $\gamma _{t}$-edge-supercritical when $\gamma_{t}(G+e)=\gamma_{t}(G)-2$ for every $e\in E(\overline{G})\neq\emptyset$. Similarly, by \thref{tR bounds}, a graph $G$ with no isolated vertices is $\gamma _{tR}$-edge-supercritical when $\gamma_{tR}(G+e)=\gamma_{tR}(G)-2$ for every $e\in E(\overline{G})\neq\emptyset$. We begin with a result by Haynes, Mynhardt and Van der Merwe \cite{HMM} characterizing $\gamma_{t}$-edge-supercritical graphs, as well as the lemma required to prove this result.      

\begin{lemma}
\thlabel{t super lemma} \emph{\cite{HMM}}
If $G$ is a graph with no isolated vertices and $u,v\in V(G)$ such that $d(u,v)=2$, then $\gamma_{t}(G)-1\leq\gamma_{t}(G+uv)$.
\end{lemma}

\begin{prop}
\thlabel{t super} \emph{\cite{HMM}} A graph $G$ is $\gamma _{t}$-edge-supercritical if and only if $G$ is the union of two or more non-trivial complete graphs.
\end{prop}

\vspace{0mm}

Lampman et al. \cite{LMO} considered whether an analogous result holds for $\gamma _{tR}$-edge-supercritical graphs. They determined that a result analogous to \thref{t super lemma} does not hold with respect to total Roman domination, and thus, even if a result similar to \thref{t super} holds, it cannot be proved via the technique employed by Haynes et al. in \cite{HMM}. However, they did establish that an analogous sufficient condition does hold for $\gamma_{tR}$-edge-supercritical graphs, which we now present.

\begin{prop} 
\thlabel{no 5-super}
\emph{\cite{LMO}}
\vspace{-8mm} 
\setlist{nolistsep}
\begin{enumerate}
\item[]
\item[$(i)$] There are no $5$-$\gamma_{tR}$-edge-supercritical graphs.
\item[$(ii)$] If $G$ is the disjoint union of $k\geq2$ complete graphs, each of order at least $3$, then $G$ is $3k$-$\gamma_{tR}$-edge-supercritical. 
\end{enumerate}
\end{prop} 

\vspace{0mm} 

Lampman et al. \cite{LMO} left the existence of connected $\gamma_{tR}$-edge-supercritical graphs as an open problem, which we investigate here. We begin by demonstrating the existence of connected $2n$-$\gamma_{tR}$-\linebreak edge-supercritical graphs for $n\geq4$. 

\begin{prop}
\thlabel{super cor} If $G=\text{cor}(K_n)$ for $n\geq4$, then $G$ is $\gamma_{tR}$-edge-supercritical.
\end{prop}

\noindent \emph{Proof.} By \thref{tR=n}, $\gamma_{tR}(G)=2n$. Label the vertices of $G$ such that $u_1,u_2,...,u_n$ are the pendant vertices with support vertices $w_1,w_2,...,w_n$, respectively. Consider $uv\in E(\overline{G})$. Then at least one of $u$ and $v$ has degree $1$ in $G$; say $\text{deg}_G(u)=1$. Note that we may assume $u=u_1$, without loss of generality. We consider two cases:
\setlist{nolistsep}
\begin{enumerate}
	\item [Case 1:] Suppose $v=u_2$ (without loss of generality). Consider $f:V(G)\rightarrow \{0,1,2\}$ defined by $f(u_1)=f(u_2)=1$, $f(w_3)=f(w_4)=...=f(w_n)=2$, and $f(z)=0$ for all other $z\in V(G)$.
	\item [Case 2:] Suppose $v=w_2$ (without loss of generality). Consider $f:V(G)\rightarrow \{0,1,2\}$ defined by $f(w_2)=f(w_3)=...=f(w_n)=2$, and $f(z)=0$ for all other $z\in V(G)$.
\end{enumerate}
In either case, $f$ is a TRD-function on $G+uv$ with $\omega(f)=2n-2$. Hence $G$ is $\gamma_{tR}$-edge-supercritical.~$\square $ 

\vspace{0mm}

Having proved the existence of connected $\gamma_{tR}$-edge-supercritical graphs, we now present the following necessary condition for a graph $G$ to be $\gamma _{tR}$-edge-supercritical. 

\begin{prop}
\thlabel{super endpaths} If $G$ is a $\gamma_{tR}$-edge-supercritical graph, then $G$ contains no adjacent endpaths. 
\end{prop}

\noindent \emph{Proof.} Suppose for a contradiction that $G$ contains two adjacent endpaths $w,v_1,...,v_n$ and $w,u_1,...,u_m$. Since $G$ is $\gamma_{tR}$-edge-supercritical, \thref{end deg 3} implies that $n,m\geq2$. Moreover, by \thref{no long legs}, at least one of $n$ and $m$ is equal to $2$; say $n=2$. Consider $u_1v_1\in E(\overline{G})$ and a $\gamma_{tR}$-function $f$ on $G+u_1v_1$. Since $n=2$, \thref{pendant edge} implies that $f(v_1)>0$. If $f(u_1)>0$, define $f^{\prime }:V(G)\rightarrow \{0,1,2\}$ by $f^{\prime}(w)=1$ and $f^{\prime}(x)=f(x)$ for all other $x\in V(G)$. Otherwise, if $f(u_1)=0$, then by \thref{set added edge}, $f(v_1)=2$. Thus, by \thref{pendant edge}, we may assume without loss of generality that $f(v_2)=0$. Hence $f(w)>0$. Therefore, define $f^{\prime }:V(G)\rightarrow \{0,1,2\}$ by $f^{\prime}(u_1)=1$ and $f^{\prime}(x)=f(x)$ for all other $x\in V(G)$. In either case, $f^{\prime}$ if a TRD-function on $G$ with $\omega(f^{\prime})\leq\omega(f)+1$, contradicting $G$ being $\gamma_{tR}$-edge-supercritical. Therefore $G$ contains no adjacent endpaths.~$\square $   

\vspace{0mm}

As a result of \thref{super endpaths}, every $\gamma_{tR}$-edge-supercritical graph contains a cycle, as we now show. 

\begin{cor}
\thlabel{super trees} There are no $\gamma_{tR}$-edge-supercritical trees.
\end{cor}

\noindent \emph{Proof.} Suppose for a contradiction that $T$ is a $\gamma_{tR}$-edge-supercritical tree. By \thref{tR=n,n edge-crit}, $T$ cannot be a path. Therefore $T$ contains at least one branch vertex (that is, a vertex of degree $3$ or more), and hence two adjacent endpaths, contradicting \thref{super endpaths}. Therefore, there are no $\gamma_{tR}$-edge-supercritical trees.~$\square $

\section{$5$-$\gamma_{tR}$-Edge-critical graphs}
\label{Sec: tR=5}

As seen in Section \ref{Sec: Prelim}, Lampman et al. characterized connected $4$-$\gamma_{tR}$-edge-critical graphs in \cite{LMO}. There, they also provided necessary conditions for a graph $G$ to be $5$-$\gamma_{tR}$-edge-critical (see \thref{if 5tR edge-crit}). In this section, we develop a characterization of $5$-$\gamma_{tR}$-edge-critical graphs from these necessary conditions.

\begin{prop} 
\thlabel{if 5tR edge-crit} \emph{\cite{LMO}} For any graph $G$, if $G$ is $5$-$\gamma _{tR}$-edge-critical, then $G$ is either $3$-$\gamma _{t}$-edge-critical or $G=K_{2}\cup K_{n}$ for $n\geq 3$, in which case $G$ is $4$-$\gamma _{t}$-edge-supercritical.
\end{prop}

\vspace{0.6mm}

Before characterizing $5$-$\gamma_{tR}$-edge-critical graphs, we characterize the connected graphs with total Roman domination number $5$, as follows.

\begin{thm}
\thlabel{tR=5} For a connected graph $G$, $\gamma_{tR}(G)=5$ if and only if $\gamma_{t}(G)=3$ and there exist a $\gamma(G)$-set $S$ and a $\gamma_{t}(G)$-set $T$ such that $S\subset T$.
\end{thm}

\noindent \emph{Proof.} Suppose $\gamma_{tR}(G)=5$. By \thref{TRD vs D}, $\gamma(G)\leq2$. Furthermore, by \thref{tR=3}, $G$ has no universal vertex. Therefore $\gamma(G)>1$, and thus $\gamma(G)=2$. Moreover, \thref{tR>t+1} implies that $\gamma_{t}(G)\leq3$. By \thref{tR=34}, $\gamma_{t}(G)\neq2$, and thus $\gamma_{t}(G)=3$. Now, consider a $\gamma_{tR}(G)$-function $f$ such that $|V_f^+|$ contains the minimum number of components. If $|V_f^2|=0$, then by \thref{tR=n}, $G\cong P_5$ or $G\cong C_5$. In either case, there exist a $\gamma(G)$-set $S$ and a $\gamma_{t}(G)$-set $T$ such that $S\subset T$. If $|V_f^2|=2$, then $V_f^2$ is a $\gamma(G)$-set and $V_f^+$ is a $\gamma_{t}(G)$-set, where $V_f^2\subset V_f^+$ as required. Otherwise, assume $|V_f^2|=1$; say $f(u)=2$. Since $f$ was chosen such that $|V_f^+|$ contains the minimum number of components, it is easy to see that $G[V_f^+]$ is connected. Therefore, $G[V_f^1]\cong P_3= v,w,x$ such that $uv\in E(G)$ but $uw,ux\notin E(G)$. Taking $S=\{u,w\}$ and $T=\{u,v,w\}$ gives the required result.

Conversely, suppose $\gamma_{t}(G)=3$. Then, since $S\subset T$, we have $\gamma(G)<3$. Hence $\gamma(G)=2$, as $G$ clearly has no universal vertex. Therefore, by \thref{TRD vs D}, $4\leq\gamma_{tR}(G)\leq6$. Furthermore, \thref{tR=34} implies that $\gamma_{tR}(G)\neq4$. Hence $\gamma_{tR}(G)\in\{5,6\}$. Suppose for a contradiction that $\gamma_{tR}(G)=6$, and consider a $\gamma(G)$-set $S$ and a $\gamma_{t}(G)$-set $T$ such that $S\subset T$. Since $\gamma_{t}(G)=3$, $G[T]\cong K_3$ or $G[T]\cong P_3$. Clearly $G[T]\ncong K_3$, otherwise $G[S]$ would be connected, contradicting $\gamma_{t}(G)=3$. Thus $G[T]\cong P_3$; say $G[T]$ is the path $u,v,w$. Since $S\subset T$, clearly $S=\{u,w\}$. However, the function $f:V(G)\rightarrow \{0,1,2\}$ defined by $f(u)=f(w)=2$, $f(v)=1$, and $f(y)=0$ for all other $y\in V(G)$ is then a TRD-function on $G$ with $\omega(f)=5$, contradicting $\gamma_{tR}(G)=6$. Therefore $\gamma_{tR}(G)=5$.~$\square $

\vspace{0mm}

The characterization of $5$-$\gamma_{tR}$-edge-critical graphs follows.

\begin{prop}
\thlabel{5tR edge-crit} A graph $G$ is $5$-$\gamma _{tR}$-edge-critical if and only if either $G$ is $3$-$\gamma _{t}$-edge-critical and there exist a $\gamma(G)$-set $S$ and a $\gamma_{t}(G)$-set $T$ such that $S\subset T$, or $G=K_{2}\cup K_{n}$ for $n\geq 3$, in which case $G$ is $4$-$\gamma_{t}$-edge-supercritical.
\end{prop}

\noindent \emph{Proof.} If $G$ is $5$-$\gamma _{tR}$-edge-critical, then the result follows directly from \thref{if 5tR edge-crit,tR=5}. Conversely, suppose $G$ is $3$-$\gamma_{t}$-edge-critical and there exists a $\gamma(G)$-set $S$ and a $\gamma_{t}(G)$-set $T$ such that $S\subset T$. Then $\gamma_{t}(G+e)=2$ for every $e\in E(\overline{G})$. Therefore \thref{tR=34} implies that $\gamma_{tR}(G+e)\in\{3,4\}$ for every $e\in E(\overline{G})$. Since $\gamma_{t}(G)=3$ and there exist a $\gamma(G)$-set $S$ and a $\gamma_{t}(G)$-set $T$ such that $S\subset T$, \thref{tR=5} implies that $\gamma_{tR}(G)=5$, and thus $G$ is $5$-$\gamma_{tR}$-edge-critical. Otherwise, if $G=K_{2}\cup K_{n}$ for $n\geq 3$, then $G$ is clearly $5$-$\gamma _{tR}$-edge-critical.~$\square $

\section{$6$-$\gamma_{tR}$-Edge-supercritical graphs}
\label{Sec: 6-super}

We now consider $\gamma_{tR}$-edge-supercritical graphs with total Roman domination number $6$, which, by \thref{no 5-super}, is the smallest TRD-number possible for a $\gamma_{tR}$-edge-supercritical graph. We begin by characterizing the disconnected $6$-$\gamma_{tR}$-edge-supercritical graphs.

\begin{prop}
\thlabel{disconnect 6-super} A disconnected graph $G$ is $6$-$\gamma_{tR}$-edge-supercritical if and only if $G\cong K_n\cup K_m$, where $n,m\geq3$. 
\end{prop}

\noindent \emph{Proof.} First, suppose $G$ is $6$-$\gamma_{tR}$-edge-supercritical. Since $\gamma_{tR}(H)\geq2$ for any graph $H$ without isolated vertices, with equality if and only if $H=K_{2}$, $G$ has two or three components. If $G$ has three components, then $G=K_{2}\cup K_{2}\cup K_{2}$ and $\gamma_{tR}(G+e)=6$ for any $e\in E(\overline{G})$, contradicting $G$ being $6$-$\gamma_{tR}$-edge-supercritical. Thus $G$ has two components; say $H_1$ and $H_2$. Now, either (say) $H_{1}=K_{2}$ and $\gamma_{tR}(H_{2})=4$, or $\gamma_{tR}(H_{1})=\gamma_{tR}(H_{2})=3$. In the former case, \thref{tR=3} implies that $H_{2}$ is not complete. Thus $\gamma_{tR}(H_{2}+e)\geq3$ for any edge $e\in E(\overline{H_{2}})\neq\emptyset$, contradicting our assumption that $G$ is $6$-$\gamma_{tR}$-edge-supercritical. In the latter case, $H_{i}$ has a universal vertex for $i=1,2$. If $H_{i}$ is not complete, then $\gamma_{tR}(H_{i}+e)=3$, and thus $\gamma_{tR}(G+e)=6$, for each edge $e\in E(\overline{H_{i}})\neq\emptyset$. We conclude that $H_{1}$ and $H_{2}$ are complete graphs of order at least $3$, as required. The converse follows directly from \thref{no 5-super}.~$\square $ 
        
\vspace{0mm}

We now consider connected $6$-$\gamma_{tR}$-edge-supercritical graphs, beginning with a result bounding the diameter of such a graph.

\begin{prop}
\thlabel{diam 6-super} If $G$ is a connected $6$-$\gamma_{tR}$-edge-supercritical graph, then $2\leq\text{\emph{diam}}(G)\leq3$. 
\end{prop}

\noindent \emph{Proof.} Clearly $2\leq\text{diam}(G)$, otherwise $E(\overline{G})=\emptyset$ and hence $G$ is not $\gamma_{tR}$-edge-critical. Now, suppose for a contradiction that $\text{diam}(G)\geq4$. Let $u$ and $v$ be vertices such that $d(u,v)=4$; say $u,x,y,z,v$ is a $u-v$ path. Since $G$ is $6$-$\gamma_{tR}$-edge-supercritical, $\gamma_{tR}(G+uv)=4$. Consider a $\gamma_{tR}$-function $f$ on $G+uv$. By \thref{set added edge}, $\{f(u),f(v)\}\in\{\{2,2\},\{2,1\},\{2,0\},\{1,1\}\}$. If $f(u)=f(v)=1$, then, in order to totally Roman dominate $\{x,y,z\}$, there exists some vertex $w\in N_G(u)$ (without loss of generality) such that $w\in N_G(x)\cap N_G(y)\cap N_G(z)$. But then $u,w,z,v$ is a shorter $u-v$ path, a contradiction. Otherwise, if $f(u)=2$ (without loss of generality), then, in order to totally Roman dominate $\{y,z\}$, there exists some vertex $w\in N_G(u)$ such that $w\in N_G(y)\cap N_G(z)$. Again, $u,w,z,v$ is a shorter $u-v$ path, a contradiction. Therefore $\text{diam}(G)\leq3$.~$\square $

\vspace{0mm}

In Section \ref{Sec: Super}, we demonstrated the existence of connected $2n$-$\gamma_{tR}$-edge-supercritical graphs for each $n\geq4$. We now demonstrate the existence of an infinite class of $6$-$\gamma_{tR}$-edge-supercritical graphs. We define the graph $G_r$ below, and show that $G_r$ is such a graph for each $r\geq2$. Note that $\text{diam}(G_r)=3$. 

Let $G_r$ be the graph constructed from $K_{2r}$ as follows: Label the vertices of $K_{2r}$ as $u_1,u_2,...,u_r,w_1,\linebreak w_2,...,w_r$, and remove from $K_{2r}$ the perfect matching $u_iw_i$ where $1\leq i\leq r$. Add a vertex disjoint $K_3$ component to $K_{2r}$, and label the added vertices $x,y,z$. Let $z$ be adjacent to both $u_i$ and $w_i$, and $y$ be adjacent to $u_i$, for $1\leq i\leq r$. Finally, add two more vertices $u_0$ and $w_0$, such that $u_0x,u_0u_i,w_0y,w_0w_i\in E(G_r)$ for $1\leq i\leq r$. See Figure \ref{G_r}.

{\begin{figure}
	\centering
	\begin{tikzpicture}[scale=0.9]		 
		\node [std] (z) at (0,0)[label=left: $z$] {};
		\node [std] (y) at (-1,-2)[label=left: $y$] {};
		\node [std] (x) at (-1,2)[label=left: $x$] {};
		
		\node [std] (u0) at (5.5,3.5)[label=right: $u_0$] {};
		\node [std] (u1) at (2,1.5)[label=above: $u_{1}$] {};
		\node [std] (u2) at (4,1.5)[label=above: $u_{2}$] {};
		\node [std] (ur) at (10,1.5)[label=above: $u_{r}$] {};

		\node [std] (w0) at (5.5,-3.5)[label=right: $w_{0}$] {};
		\node [std] (w1) at (2,-1.5)[label=below: $w_{1}$] {};
		\node [std] (w2) at (4,-1.5)[label=below: $w_{2}$] {};
		\node [std] (wr) at (10,-1.5)[label=below: $w_{r}$] {};

		\foreach \i in {1,2}
			{\draw (u\i)--(u0);
			\draw (w\i)--(w0);
			\draw (u\i)--(z);
			\draw (w\i)--(z);
			\draw (u\i)--(y);
			\draw (u\i)--(wr);
			\draw (w\i)--(ur);}
			
		\draw (ur)--(u0);
		\draw (wr)--(w0);
		\draw (ur)--(z);
		\draw (wr)--(z);
		\draw (ur)--(y);
		
		\draw (x)--(y);
		\draw(y)--(z);
		\draw(x)--(z);
		\draw(x)--(u0);
		\draw(y)--(w0);
		\draw(u1)--(w2);
		\draw(u2)--(w1);
		\draw(w1)--(w2);
		\draw(u2)--(u1);
		
		\draw (u2)--(ur);
		\draw (w2)--(wr);
		
		\node[draw=Plum, thick, rounded corners, dashed, fit=(w1) (wr)] (W){}; 
		\node[draw=Plum, thick, rounded corners, dashed, fit=(u1) (ur)] (U){}; 
		\node[right=.1cm of W][color=Plum](W1){$K_r$};
		\node[right=.1cm of U][color=Plum](U1){$K_r$};
		
	\end{tikzpicture}	
	\caption{The graph $G_r$, where $r\geq2$.}
	\label{G_r}
\end{figure}
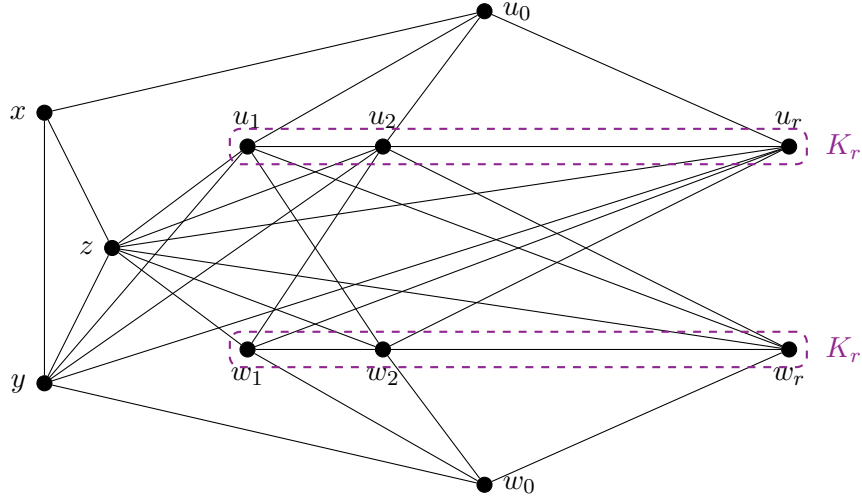}

\begin{thm}
\thlabel{G_r super} If $r\geq2$, then $G_r$ is $6$-$\gamma_{tR}$-edge-supercritical. 
\end{thm} 

\pagebreak

\noindent \emph{Proof.} We first show that $\gamma_{tR}(G_r)=6$ for $r\geq2$. By inspection, $\gamma(G)>2$. Therefore, since $\{x,y,z\}$ dominates $G_r$, $\gamma(G)=3$. Furthermore, this is a TD-set on $G$, and thus $\gamma_{t}(G)=3$. By \thref{TRD vs D}, $\gamma_{tR}(G)\leq6$. Moreover, \thref{tR=34,tR=5} imply that $\gamma_{tR}(G)>5$, and hence $\gamma_{tR}(G)=6$. 

Now, consider any edge $vv^{\prime}\in E(\overline{G_r})$. Consider the following cases:
\setlist{nolistsep}
\begin{enumerate}
	\item [Case 1:] Let $v=u_0$. Then, without loss of generality, $v^{\prime}\in\{y,z,w_0,w_1\}$. If $v^{\prime}\in\{y,w_1\}$, consider the function $f:V(G_r)\rightarrow\{0,1,2\}$ defined by $f(v^{\prime})=f(z)=2$ and $f(b)=0$ for all other $b\in V(G_r)$. Otherwise, if $v^{\prime}\in\{z,w_0\}$, consider the function $f:V(G_r)\rightarrow\{0,1,2\}$ defined by $f(v^{\prime})=f(y)=2$ and $f(b)=0$ for all other $b\in V(G_r)$.
	\item [Case 2:] Let $v=z$. Then $v^{\prime}=w_0$. Consider the function $f:V(G_r)\rightarrow\{0,1,2\}$ defined by $f(u_1)=f(z)=2$ and $f(b)=0$ for all other $b\in V(G_r)$.
	\item [Case 3:] Let $v=y$. Then, without loss of generality, $v^{\prime}=w_1$. Consider the function $f:V(G_r)\rightarrow\{0,1,2\}$ defined by $f(y)=f(u_1)=2$ and $f(b)=0$ for all other $b\in V(G_r)$.
	\item [Case 4:] Let $v=w_0$. Then, without loss of generality, $v^{\prime}\in\{x,u_1\}$. Consider the function $f:V(G_r)\rightarrow\{0,1,2\}$ defined by $f(v^{\prime})=f(z)=2$ and $f(b)=0$ for all other $b\in V(G_r)$.   
	\item [Case 5:] Let $v=x$. Then, without loss of generality, $v^{\prime}\in\{u_1,w_2\}$. Consider the function $f:V(G_r)\rightarrow\{0,1,2\}$ defined by $f(u_1)=f(w_2)=2$ and $f(b)=0$ for all other $b\in V(G_r)$.
	\item [Case 6:] Let $v=u_1$ (without loss of generality). Then $v^{\prime}=w_1$. Consider the function $f:V(G_r)\rightarrow\{0,1,2\}$ defined by $f(y)=f(u_1)=2$ and $f(b)=0$ for all other $b\in V(G_r)$.
\end{enumerate}
In each case, $f$ is a TRD-function on $G_r+vv^{\prime}$ with $\omega(f)=4$. Therefore, by \thref{tR bounds}, $\gamma_{tR}(G_r+e)=4$ for any $e\in E(\overline{G_r})$. Thus $G_r$ is $6$-$\gamma_{tR}$-edge-supercritical.~$\square $

\begin{cor}
For $r\geq2$, there exists a connected $6$-$\gamma_{tR}$-edge-supercritical graph on $5+2r$ vertices. 
\end{cor}

\section{$\gamma_{tR}$-Edge-removal-critical graphs}
\label{Sec: ER-critical}

We now consider the effect that the removal of an edge has on the total Roman domination number of a graph. The following observations follow directly from \thref{tR bounds,set added edge}, and \thref{pendant edge}.

\begin{obs}
\thlabel{tR ER bounds} 
Given a graph $G$ with no isolated vertices, if ${uv\in E_P(G)}$, then $\gamma
_{tR}(G)\leq \gamma _{tR}(G-uv)\leq \gamma _{tR}(G)+2$.
\end{obs}

\begin{obs}
\thlabel{ER set} For a graph $G$ with no isolated vertices, if $uv\in E(G)$ is a TRD-ER-critical edge, then, for any $\gamma_{tR}(G)$-function $f$, $\{f(u),f(v)\}\in\{\{0,2\},\{1,2\},\{2,2\},\{1,1\}\}$.
\end{obs}

%\vspace{0mm}
%\pagebreak

As with TRD-ER-critical edges, we now present a result restricting the possible values assigned to the vertices of a TRD-ER-supercritical edge $e\in E(G)$ by a $\gamma_{tR}$-function $f$ on $G$.

\begin{prop}
\thlabel{super ER set} For a graph $G$ with no isolated vertices, if $uv\in E(G)$ is a TRD-ER-supercritical edge, then there exists a $\gamma _{tR}(G)$-function $f$ such that $\{f(u),f(v)\}\in\{\{2,2\},\{2,0\},\{1,1\}\}$.
\end{prop}

\noindent \emph{Proof.} Let $G^{\prime}=G-uv$. By \thref{ER set}, $\{f(u),f(v)\}\in\{\{2,2\},\{2,1\},\{2,0\},\{1,1\}\}$ for any $\gamma _{tR}(G)$-function $f$. Suppose for a contradiction that $\{f(u),f(v)\}=\{1,2\}$ for every $\gamma _{tR}(G)$-function $f$, and consider one such function. Say $f(u)=2$ and $f(v)=1$. Then by \thref{pendant edge}, $\text{deg}(u)>1$ and $\text{deg}(v)>1$. Now, $f$ is a RD-function on $G^{\prime}$, with $u$ and $v$ being the only possible isolated vertices in $G^{\prime}[V_f^+]$. Note that at least one of $u$ and $v$ must be isolated in $G^{\prime}[V_f^+]$, otherwise $f$ is also a TRD-function on $G^{\prime}$, contradicting $uv$ being TRD-ER-critical. 

Suppose for a contradiction that $v$ is isolated in $G^{\prime}[V_{f}^{+}]$. That is, $f(x)=0$ for all $x\in N_G(v)-\{u\}$. Since $\text{deg}(u)>1$, there exists some $w\in N_G(u)-\{v\}$. But $f(w)=0$ for each such $w$, otherwise $f^{\prime}:V(G)\rightarrow \{0,1,2\}$ defined by $f^{\prime}(v)=0$ and $f^{\prime}(z)=f(z)$ for all other $z\in V(G)$ would be a TRD-function on $G$, contradicting the minimality of $f$. That is, $u$ is also isolated in $G^{\prime}[V_{f}^{+}]$. But then $g:V(G)\rightarrow \{0,1,2\}$ defined by $g(v)=0$, $g(w)=1$ for some $w\in N(u)-\{v\}$, and $g(z)=f(z)$ for all other $z\in V(G)$ is a $\gamma_{tR}(G)$-function with $g(u)=2$ and $g(v)=0$, contradicting our assumption. 

Therefore $u$ is the only isolated vertex in $G^{\prime}[V_{f}^{+}]$. But then $g:V(G)\rightarrow \{0,1,2\}$ defined by $g(w)=1$ for some $w\in N_G(u)-\{v\}$ and $g(z)=f(z)$ for all other $z\in V(G)$ is a TRD-function on $G$ with $\omega(g)=\omega(f)+1$, contradicting $uv$ being a TRD-ER-supercritical edge.~$\square$

\begin{cor}
\thlabel{super set} For a graph $G$ with no isolated vertices, if $uv\in E(\overline{G})$ is a TRD-supercritical edge, then there exists a $\gamma _{tR}(G+uv)$-function $f$ such that $\{f(u),f(v)\}\in\{\{2,2\},\{2,0\},\{1,1\}\}$.
\end{cor} 

\vspace{0mm}

We now consider $\gamma _{t}$-ER-critical and $\gamma_{tR}$-ER-critical graphs. Recall that a graph $G$ with no isolated vertices is $\gamma _{t}$-ER-critical if $\gamma_{t}(G+e)>\gamma_{t}(G)$ for every $e\in E(G)$, and similarly $\gamma _{tR}$-ER-critical if $\gamma_{tR}(G+e)>\gamma_{tR}(G)$ for every $e\in E(G)$. Connected $\gamma_{t}$-ER-critical graphs $G$ were characterized in \cite{DHH}. There, Desormeaux et al. defined $\mathcal{T}$ to be the family of trees $T$ such that $T$ is either a nontrivial star, or a double star, or can be obtained from a subdivided star by adding zero or more pendant edges to the non-leaf vertices.

\begin{prop}
\thlabel{TD-ER-crit} \emph{\cite{DHH}}
A connected graph $G$ is $\gamma_{t}$-ER-critical if and only if $G\in \mathcal{T}$. 
\end{prop}

\vspace{0mm}

Note that a disconnected graph $G$ is $\gamma_{t}$-ER-critical if and only if each component of $G$ is itself $\gamma_{t}$-ER-critical. As a result, \thref{TD-ER-crit} provides the following characterization of all $\gamma _{t}$-ER-critical graphs.

\begin{obs}
\thlabel{TD-ER-crit disc} 
A graph $G$ is $\gamma_{t}$-ER-critical if and only if $G$ is the union of $k\geq1$ graphs $G_i\in \mathcal{T}$, for $1\leq i\leq k$. 
\end{obs}

\vspace{0mm}

We investigate whether a similar characterization holds for $\gamma_{tR}$-ER-critical graphs. Note that as with $\gamma_{t}$-ER-critical graphs, a disconnected graph $G$ is $\gamma_{tR}$-ER-critical if and only if each component of $G$ is itself $\gamma_{tR}$-ER-critical. Similarly, a disconnected graph $G$ is $\gamma_{tR}$-ER-supercritical if and only if each component of $G$ is itself $\gamma_{tR}$-ER-supercritical. As a result, we focus specifically on connected $\gamma_{tR}$-ER-critical and $\gamma_{tR}$-ER-supercritical graphs. We begin with a result restricting the values that a $\gamma_{tR}(G)$-function $f$ can assign to the vertices of a $\gamma_{tR}$-ER-critical graph based on their degree.

\begin{prop}
\thlabel{ER crit min deg} Let $G$ be a connected $\gamma_{tR}$-ER-critical graph. For any $\gamma_{tR}$-function $f$ on $G$, if $f(u)=0$, then $\text{\emph{deg}}(u)=1$. Moreover, $\delta(G)=1$.
\end{prop}

\noindent \emph{Proof.} Let $G$ be a connected $k$-$\gamma_{tR}$-ER-critical graph of order $n$, and $f$ a $\gamma_{tR}(G)$-function. Suppose for a contradiction that there exists $u\in V(G)$ such that $f(u)=0$ and $\text{deg}(u)\geq2$. Then there exist $v,w\in N_G(u)$. By \thref{ER set}, $f(v)=f(w)=2$. But then $f$ is also a TRD-function on $G-uv$, contradicting $uv$ being TRD-ER-critical. Hence $\text{deg}(u)=1$. Now, if $\delta(G)\geq2$, then $V_f^1=V(G)$; that is, $k=n$. But then \thref{tR ER bounds} implies that $\gamma_{tR}(G-e)=n=k$ for all $e\in E(G)$, contradicting our assumption that $G$ is $\gamma_{tR}$-ER-critical. Hence $\delta(G)=1$.~$\square$ 

\vspace{0mm} 

Note that \thref{ER crit min deg} implies that every component of a $\gamma_{tR}$-ER-critical graph contains at least one degree $1$ vertex. We now present a result demonstrating that a connected $\gamma_{tR}$-ER-critical graph $G$ cannot contain any cycles.

\begin{prop} 
\thlabel{ER crit trees}
If $G$ is a connected $\gamma_{tR}$-ER-critical graph, then $G$ is a tree. 
\end{prop}

\noindent \emph{Proof.} Suppose for a contradiction that $G$ is a connected $\gamma_{tR}$-ER-critical graph which contains a cycle; say $v_1,v_2,...,v_k,v_1$, for $k\geq3$. Consider a $\gamma_{tR}$-function $f$ on $G$. By \thref{ER crit min deg}, $f(v_i)>0$ for $1\leq i\leq k$. But then $f$ is also a TRD-function on $G-v_1v_2$, contradicting $G$ being $\gamma_{tR}$-ER-critical. Hence $G$ cannot contain a cycle, and thus, since $G$ is connected, $G$ is a tree.~$\square$  

\vspace{0mm}

Our next result restricts the distance between any two vertices of a $\gamma_{tR}$-ER-critical graph $G$ which are in $V_f^+$ for some $\gamma_{tR}(G)$-function $f$. 

\begin{prop} 
\thlabel{ER crit max dist}
Let $G$ be a connected $\gamma_{tR}$-ER-critical graph. If $u,v\in V(G)$ and $f$ is a $\gamma_{tR}$-function on $G$ such that $f(u)>0$ and $f(v)>0$, then $d(u,v)\leq 2$. 
\end{prop}

\noindent \emph{Proof.} Let $G$ be a connected $\gamma_{tR}$-ER-critical graph. Then, by \thref{ER crit trees}, $G$ is a tree. Let $f(u)>0$ and $f(v)>0$, and suppose for a contradiction that $u,w_1,...,w_k,v$ is the unique path from $u$ to $v$, where $k\geq2$. Consider a $\gamma_{tR}$-function $f$ on $G$. Then \thref{ER crit min deg} implies that $f(w_i)>0$ for all $1\leq i\leq k$. But then $f$ is a TRD-function on $G-w_1w_2$, contradicting $G$ being $\gamma_{tR}$-ER-critical.~$\square$ 

\begin{cor}
\thlabel{ER crit max diam}
Let $G$ be a connected $\gamma_{tR}$-ER-critical graph. If $u,v\in V(G)$ such that $\text{\emph{deg}}(u)>1$ and $\text{\emph{deg}}(v)>1$, then $d(u,v)\leq2$. Moreover, $\text{\emph{diam}}(G)\leq4$.
\end{cor}

\vspace{0mm}

We now present a characterization of the graphs $G$ which are $\gamma_{tR}$-edge-removal-critical. Consider for a moment a star graph $S_n$, which is defined to be the complete bipartite graph $K_{1,n}$, with $n\geq1$. Let $\mathcal{F}_n$ be the family of graphs constructed from $S_n$ by appending $k_1,k_2,...,k_n$ (where $k_1\geq k_2\geq ... \geq k_n\geq 0$) pendant vertices to each pendant vertex of $S_n$. In what follows, we label the vertices of a graph $G\in\mathcal{F}_n$ as follows: Let $c$ be the central vertex (choosing either central vertex in the case of $S_{1}$), and $u_i$ ($1\leq i\leq n$) the pendant vertices, in the original star $S_n$. For each such vertex $u_i$, let $v_{i,1},v_{i,2},...,v_{i,k_i}$ be the pendant vertices added to $u_i$. See Figure \ref{F_n}.

{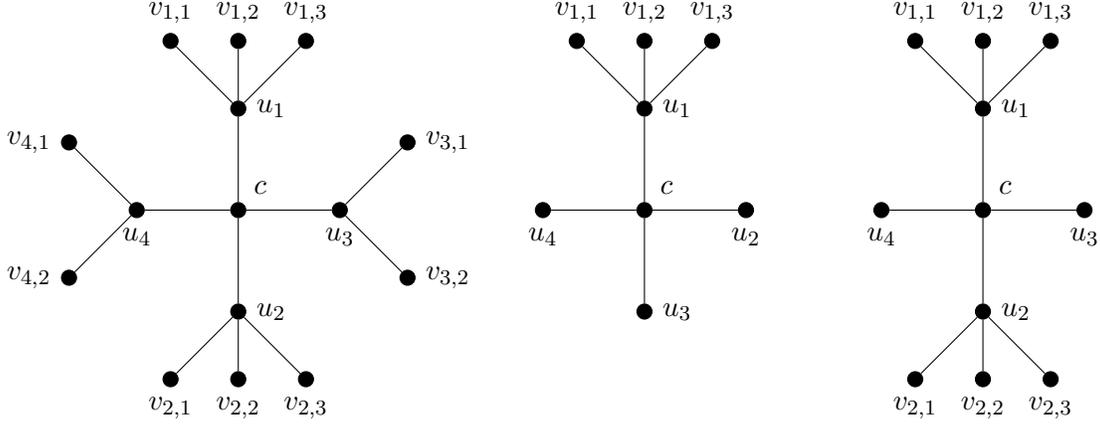
\begin{figure}
	\centering
	\begin{tikzpicture}[scale=0.9]
		\begin{scope}[shift={(-9,0)}]	 
		\node [std] (c) at (0,0)[label=above right: $c$] {};
		
		\node [std] (u1) at (0,1.5)[label=right: $u_1$] {};
		\node [std] (u2) at (0,-1.5)[label=right: $u_2$] {};
		\node [std] (u3) at (1.5,0)[label=below: $u_3$] {};
		\node [std] (u4) at (-1.5,0)[label=below: $u_4$] {};
		
		\node [std] (v11) at (-1,2.5)[label=above: $v_{1,1}$] {};
		\node [std] (v12) at (0,2.5)[label=above: $v_{1,2}$] {};
		\node [std] (v13) at (1,2.5)[label=above: $v_{1,3}$] {};
				
		\node [std] (v31) at (2.5,1)[label=right: $v_{3,1}$] {};
		\node [std] (v32) at (2.5,-1)[label=right: $v_{3,2}$] {};
			
		\node [std] (v21) at (-1,-2.5)[label=below: $v_{2,1}$] {};
		\node [std] (v22) at (0,-2.5)[label=below: $v_{2,2}$] {};
		\node [std] (v23) at (1,-2.5)[label=below: $v_{2,3}$] {};
		
		\node [std] (v41) at (-2.5,1)[label=left: $v_{4,1}$] {};
		\node [std] (v42) at (-2.5,-1)[label=left: $v_{4,2}$] {};

		\foreach \i in {1,2,3,4}
			{\draw (c)--(u\i);
			\draw (u\i)--(v\i1);
			\draw (u\i)--(v\i2);}
			
		\draw (u1)--(v13);
		\draw (u2)--(v23);
		\end{scope}

		\begin{scope}[shift={(-3,0)}]
		\node [std] (c') at (0,0)[label=above right: $c$] {};
		
		\node [std] (u1') at (0,1.5)[label=right: $u_1$] {};
		\node [std] (u2') at (1.5,0)[label=below: $u_2$] {};
		\node [std] (u3') at (0,-1.5)[label=right: $u_3$] {};
		\node [std] (u4') at (-1.5,0)[label=below: $u_4$] {};
				
		\node [std] (v11') at (-1,2.5)[label=above: $v_{1,1}$] {};
		\node [std] (v12') at (0,2.5)[label=above: $v_{1,2}$] {};
		\node [std] (v13') at (1,2.5)[label=above: $v_{1,3}$] {};

		\foreach \i in {1,2,3,4}
			{\draw (c')--(u\i');}
			
		\draw (u1')--(v11');
		\draw (u1')--(v12');
		\draw (u1')--(v13');
		\end{scope}

		\begin{scope}[shift={(2,0)}]	 
		\node [std] (c'') at (0,0)[label=above right: $c$] {};
		
		\node [std] (u1'') at (0,1.5)[label=right: $u_1$] {};
		\node [std] (u2'') at (0,-1.5)[label=right: $u_2$] {};
		\node [std] (u3'') at (1.5,0)[label=below: $u_3$] {};
		\node [std] (u4'') at (-1.5,0)[label=below: $u_4$] {};
		
		\node [std] (v11'') at (-1,2.5)[label=above: $v_{1,1}$] {};
		\node [std] (v12'') at (0,2.5)[label=above: $v_{1,2}$] {};
		\node [std] (v13'') at (1,2.5)[label=above: $v_{1,3}$] {};
		
		\node [std] (v21'') at (-1,-2.5)[label=below: $v_{2,1}$] {};
		\node [std] (v22'') at (0,-2.5)[label=below: $v_{2,2}$] {};
		\node [std] (v23'') at (1,-2.5)[label=below: $v_{2,3}$] {};

		\foreach \i in {1,2}
			{\draw (c'')--(u\i'');
			\draw (u\i'')--(v\i1'');
			\draw (u\i'')--(v\i2'');
			\draw (u\i'')--(v\i3'');}
			
		\draw (u3'')--(c'');
		\draw (u4'')--(c'');
		\end{scope}		
	\end{tikzpicture}	
	\caption{Examples of graphs in $\mathcal{F}_4$}
	\label{F_n}
\end{figure}}

\begin{thm}
\thlabel{ER crit connect} A connected graph $G$ with no isolated vertices is $\gamma_{tR}$-ER-critical if and only if $G$ is a member of $\mathcal{F}_n$, for some $n\geq1$, such that $k_1,k_2,...,k_n\neq1$.  
\end{thm}

\noindent \emph{Proof.} Let $G$ be a connected $\gamma_{tR}$-ER-critical graph. We begin by showing that $G\in\mathcal{F}_n$ for $n\geq1$. By \thref{ER crit trees}, $G$ is a tree. Let $S=\{v\in V(G):\text{deg}_G(v)>1\}$. If $G\cong S_n$ for $n\geq1$, then $G\in\mathcal{F}_n$ as required. So assume $|S|\geq2$. We claim that $G[S]\cong S_n$ for $n\geq1$. Suppose for a contradiction that $E_P(G[S])\neq\emptyset$. Then there exist $u,v\in S$ such that $d(u,v)\geq3$. But then, by definition of $S$, $\text{diam}(G)>4$, contradicting \thref{ER crit max diam}. Hence $G[S]\cong S_n$ for $n\geq1$, and thus $G\in\mathcal{F}_n$. 

Now, consider a graph $G\in \mathcal{F}_n$ for some $n\geq1$. In what follows, let the vertices of $G$ be labelled as described in the definition of $\mathcal{F}_n$.  

\begin{enumerate}
	\item [Case 1:] Suppose $G\in \mathcal{F}_n$ for some $n\geq1$ such that $k_1,k_2,...,k_n\neq1$. If $G$ is a star or a double star, then $G$ is clearly $\gamma_{tR}$-ER-critical. Therefore, assume $n\geq2$ and $k_1\geq k_2\geq2$. Let $2\leq l\leq n$ be such that $k_i=0$ if and only if $i>l$. Note that $E_P(G)=\{cu_i: 1\leq i\leq l\}$. We consider two cases.
	\vspace{2mm}
	\begin{enumerate}
		\item [Case 1A:] Suppose $l=n$. Then it can be easily seen that $f:V(G)\rightarrow\{0,1,2\}$ defined by $f(c)=1$, $f(u_i)=2$ for all $1\leq i\leq n$, and $f(b)=0$ for all other $b\in V(G)$ is a $\gamma_{tR}(G)$-function. If $n\geq3$, then $G-cu_i$ ($1\leq i\leq n$) is the disjoint union of a star on at least $3$ vertices with a graph $H\in \mathcal{F}_{n-1}$, where $n-1\geq 2$. Otherwise, if $n=2$, $G-cu_i$ is the disjoint union of two stars, each on at least $3$ vertices. In either case, it can be easily seen that $f^{\prime}:V(G)\rightarrow\{0,1,2\}$ defined by $f^{\prime}(v_{i,1})=1$ and $f^{\prime}(z)=f(z)$ for all other $z\in V(G)$ is a $\gamma_{tR}(G-cu_i)$-function with $\omega(f^{\prime})=\omega(f)+1$, for each $1\leq i\leq n$.
		\vspace{1mm}
		\item [Case 1B:] Suppose $l<n$. Then it can be easily seen that $f:V(G)\rightarrow\{0,1,2\}$ defined by $f(c)=2$, $f(u_i)=2$ for all $1\leq i\leq l$, and $f(b)=0$ for all other $b\in V(G)$ is a $\gamma_{tR}(G)$-function. Since $2\leq l<n$, we have $n\geq3$. Hence $G-cu_i$ ($1\leq i\leq l$) is the disjoint union of a star on at least $3$ vertices with a graph $H\in \mathcal{F}_{n-1}$, where $n-1\geq 2$. Thus, it can be easily seen that $f^{\prime}:V(G)\rightarrow\{0,1,2\}$ defined by $f^{\prime}(v_{i,1})=1$ and $f^{\prime}(z)=f(z)$ for all other $z\in V(G)$ is a $\gamma_{tR}(G-cu_i)$-function with $\omega(f^{\prime})=\omega(f)+1$, for each $1\leq i\leq l$. 
\end{enumerate}

\vspace{2mm}

Therefore, in each case, $G$ is $\gamma_{tR}$-ER-critical, as required.

\vspace{2mm}

	\item [Case 2:] Otherwise, suppose $G\notin\mathcal{F}_n$ for any $n\geq1$ such that $k_1,k_2,...,k_n\neq1$. Thus $G\in\mathcal{F}_n$ for $n\geq1$ where $k_i=1$ for some $1\leq i\leq n$. If $n=1$, then $G$ is also a member of $\mathcal{F}_{2}$. Therefore, it suffices to consider $n\geq2$. Consider a $\gamma_{tR}(G)$-function $f$ such that $|V_f^2|$ is a minimum. Then $f(u_i)=f(v_{i,1})=1$. Moreover, \thref{ER crit min deg} implies that $f(c)>0$. Suppose first that $n=2$, and let $j\neq i$. If $k_j=0$, then $f(u_j)=f(c)=1$ by our choice of $f$. If $k_j\geq1$, then $f(u_j)>0$ by \thref{ER crit min deg}. Otherwise, suppose $n\geq3$. If $k_j=0$ for all $j\neq i$, then $G$ is also a member of $\mathcal{F}_2$ with $k_1=n-1\geq2$ and $k_2=0$, contradicting our assumption. Hence there exists $j\neq i$ such that $k_j\geq1$, and thus by \thref{ER crit min deg}, $f(u_j)>0$. Note that, in each case, there exists $j\neq i$ such that $f(u_j)>0$. But $u_j,c,u_i,v_{i,1}$ is a path in $G$, contradicting \thref{ER crit max dist}. Hence $G$ is not  $\gamma_{tR}$-ER-critical.~$\square$ 
\end{enumerate}

\begin{cor}
\thlabel{ER crit} A graph $G$ with no isolated vertices is $\gamma_{tR}$-ER-critical if and only if $G$ is the disjoint union of $m\geq1$ graphs $G_i\in \mathcal{F}_{n_i}$, for some $n_i\geq1$ such that $k_1,k_2,...,k_{n_i}\neq1$, for $1\leq i\leq m.$ 
\end{cor}

\section{$\gamma_{tR}$-Edge-removal-supercritical graphs}
\label{Sec: ER-supercritical}

Having classified $\gamma_{tR}$-ER-critical graphs, we now classify the graphs $G$ which are $\gamma_{tR}$-ER-supercritical.

\begin{thm}
\thlabel{ER super connect} A connected graph $G$ with no isolated vertices is $\gamma_{tR}$-ER-supercritical if and only if $G$ is either a non-trivial star, or a double star where each non-pendant vertex has degree at least $3$.
\end{thm}

\noindent \emph{Proof.} Suppose $G$ is $\gamma_{tR}$-ER-supercritical. If $E_P(G)=\emptyset$, then $G=S_n$ for $n\geq1$. Otherwise, assume $E_P(G)\neq\emptyset$. We claim that $|E_P(G)|=1$. Suppose for a contradiction that $|E_P(G)|\geq2$, and consider a path $u,v,w,x,y$ in $G$. Let $f$ be a $\gamma_{tR}(G)$-function. Then, by \thref{ER crit min deg}, $v,w,x\in V_f^+$. Moreover, since \thref{ER crit trees} implies that $G$ is a tree, by \thref{ER crit max diam}, $\text{deg}(u)=\text{deg}(y)=1$. Thus \thref{pendant edge} implies that $f(u)\leq1$ and $f(y)\leq1$. But then $g:V(G)\rightarrow\{0,1,2\}$ defined by $g(u)=1$ and $g(z)=f(z)$ for all other $z\in V(G)$ is a $\gamma_{tR}(G-vw)$-function with $\omega(g)\leq\omega(f)+1$, contradicting $vw$ being TRD-ER-supercritical. Hence $|E_P(G)|=1$, and thus $G$ is a double star.

Conversely, $G=S_n$ for $n\geq1$ is, by definition, $\gamma_{tR}$-ER-supercritical. Otherwise, suppose $G$ is a double star. Then $\gamma_{tR}(G)=4$. Moreover, $E_P(G)=\{uv\}$ where $u$ and $v$ are the two non-pendant vertices. If each non-pendant vertex has degree at least $3$, then by \thref{tR=3}, $\gamma_{tR}(G-uv)=6$, since the removal of the non-pendant edge disconnects the graph into two stars each on at least $3$ vertices. Therefore $G$ is $\gamma_{tR}$-ER-supercritical. Otherwise, if $G$ has a non-pendant vertex of degree $2$, then $\gamma_{tR}(G-uv)\leq5$, since since the removal of the non-pendant edge disconnects the graph into two stars, at least one of which is on only two vertices. Therefore $G$ is not $\gamma_{tR}$-ER-supercritical.~$\square $

\begin{cor}
\thlabel{ER super} A graph $G$ with no isolated vertices is $\gamma_{tR}$-ER-supercritical if and only if $G$ is the disjoint union of $m\geq1$ graphs $G_i$ such that, for each $1\leq i\leq m$, $G_i$ is either a non-trivial star, or a double star where each non-pendant vertex has degree at least $3$.
\end{cor}

\vspace{0mm}

We conclude this section by observing a link between $\gamma_{tR}$-ER-supercritical and $\gamma_{tR}$-edge-stable graphs, which follows directly from the previous corollary.

\begin{cor}
\thlabel{ER-super gives stable} If $G$ is a $k$-$\gamma_{tR}$-edge-removal-supercritical graph, then $G$ is $k$-$\gamma_{tR}$-edge-stable.
\end{cor}

\section{$\gamma_{tR}$-Edge-removal-stable graphs}
\label{Sec: ER-stable}

We now consider graphs $G$ which are $\gamma_{tR}$-edge-removal-stable. Recall that a graph $G$ is $\gamma_{tR}$-ER-stable when $\gamma_{tR}(G-e)=\gamma_{tR}(G)$ for all $e\in E(G)$. We begin with two observations that follow directly from the definitions of $\gamma_{tR}(G-e)$, where $e$ is a pendant edge of $G$, and a TRD-ER-stable edge, respectively.

\begin{obs}
\thlabel{ER-stable deg>1} If $G$ is a $\gamma_{tR}$-ER-stable graph, then $\delta(G)>1$.
\end{obs}

\begin{obs}
\thlabel{TRD-stable f} If $G$ is a $\gamma_{tR}$-ER-stable graph, then for any $e\in E(G)$ there exists a $\gamma_{tR}$-function $f$ on $G$ such that $f$ is also a $\gamma_{tR}(G-e)$-function.
\end{obs}

\vspace{0mm}

Consider again the graph $G_r$ defined in Section \ref{Sec: 6-super}. There, we showed that, for $r\geq2$, $G_r$ is $6$-$\gamma_{tR}$-edge-supercritical. In addition, it can be shown that $G_r$ is $\gamma_{tR}$-ER-stable. Furthermore, note that the union of $k\geq2$ complete graphs each of order at least $3$ is both $3k$-$\gamma_{tR}$-edge-supercritical (by \thref{no 5-super}) and $3k$-$\gamma_{tR}$-ER-stable (by \thref{tR=n}). Similarly, $\text{cor}(K_n)$ for $n\geq4$ is $2n$-$\gamma_{tR}$-edge-supercritical (by \thref{n edge-crit}) and every non-pendant edge $e\in E(G)$ is TRD-ER-stable (by \thref{tR=n}). In light of these results, we present the following theorem.

\begin{thm}
\thlabel{super to ER-stable} If $G$ is a $\gamma_{tR}$-edge-supercritical graph, then every non-pendant edge $e\in E(G)$ is TRD-ER-stable.
\end{thm}

\noindent \emph{Proof.} Let $G$ be a $\gamma_{tR}$-edge-supercritical graph. Then $G$ contains no $K_2$ components. Suppose for a contradiction that there exists a non-pendant edge $uw\in E(G)$ that is TRD-ER-critical. Then $\text{deg}(u)\geq2$ and $\text{deg}(w)\geq2$. Let $v\in N_G(w)-\{u\}$.

\noindent\textbf{Claim:} \emph{$N_G[u]\neq N_G[w]$.}

%\pagebreak

\emph{Proof of Claim:} Suppose for a contradiction that $N_G[u]=N_G[w]$. Let $S=N_G[u]-\{u,w\}$. Then $v\in S$. Consider a $\gamma_{tR}(G)$-function $f$. By \thref{ER set}, $\{f(u),f(w)\}\in\{\{2,2\},\{2,1\},\{2,0\},\{1,1\}\}$. We claim that $G[S]$ has no universal vertex. Suppose for a contradiction that $v$ is a universal vertex of $G[S]$. Note that possibly $S=\{v\}$. If $f(u)=f(w)=1$, consider $f^{\prime }:V(G)\rightarrow \{0,1,2\}$ defined by $f^{\prime}(u)=f^{\prime}(w)=0$, $f^{\prime}(v)=2$ and $f^{\prime}(b)=f(b)$ for all other $b\in V(G)$. Otherwise, if $f(u)=2$ (without loss of generality), consider $f^{\prime }:V(G)\rightarrow \{0,1,2\}$ defined by $f^{\prime}(u)=f(v)$, $f^{\prime}(v)=f(u)$ and $f^{\prime}(b)=f(b)$ for all other $b\in V(G)$. In either case, $f^{\prime}$ is a $\gamma_{tR}(G)$-function. Moreover, $f^{\prime}$ is also a TRD-function on $G-uw$, contradicting $uw$ being TRD-ER-critical. Therefore $G[S]$ has no universal vertex, and thus there exists some vertex $x\in S-\{v\}$ such that $vx\in E(\overline{G})$. 

Now, consider a $\gamma_{tR}$-function $g$ on $G+vx$. By \thref{set added edge}, $\{g(x),g(v)\}\in\{\{2,2\},\{2,1\},\{2,0\},\linebreak \{1,1\}\}$. If $g(x)>0$ and $g(v)>0$, then $g^{\prime }:V(G)\rightarrow \{0,1,2\}$ defined by $g^{\prime}(u)=1$ and $g^{\prime}(b)=g(b)$ for all other $b\in V(G)$ is a TRD-function on $G$ with $\omega(g^{\prime})\leq \omega(g)+1$, contradicting $vx$ being TRD-supercritical. Hence $\{g(x),g(v)\}=\{2,0\}$; say $g(v)=2$ and $g(x)=0$ (without loss of generality). Then $g(u)=g(w)=0$, otherwise $g^{\prime }:V(G)\rightarrow \{0,1,2\}$ defined by $g^{\prime}(x)=1$ and $g^{\prime}(b)=g(b)$ for all other $b\in V(G)$ would be a TRD-function on $G$ with $\omega(g^{\prime})=\omega(g)+1$, contradicting $vx$ being TRD-supercritical. Hence $h:V(G)\rightarrow \{0,1,2\}$ defined by $h(u)=h(x)=1$ and $h(b)=g(b)$ for all other $b\in V(G)$ is a $\gamma_{tR}(G)$-function, which, since $uw$ is TRD-ER-critical, contradicts \thref{ER set}. Therefore, $N_G[u]\neq N_G[w]$.~$_{(\square)}$ 

As a result of the above claim, we can choose $v\in N_G(w)-\{u\}$ such that $uv\in E(\overline{G})$. Now, consider a $\gamma_{tR}$-function $f$ on $G+uv$. By \thref{set added edge}, $\{f(u),f(v)\}\in\{\{2,2\},\{2,1\},\{2,0\},\{1,1\}\}$. If $f(u)>0$ and $f(v)>0$, then $f^{\prime }:V(G)\rightarrow \{0,1,2\}$ defined by $f^{\prime}(w)=1$ and $f^{\prime}(b)=f(b)$ for all other $b\in V(G)$ is a TRD-function on $G$ with $\omega(f^{\prime})\leq\omega(f)+1$, contradicting $uv$ being TRD-supercritical. Hence $\{f(u),f(v)\}=\{2,0\}$. We show that $f(u)=0$ and $f(v)=2$.

Suppose for a contradiction that $f(u)=2$ and $f(v)=0$. Clearly $f(w)=0$, otherwise $f^{\prime }:V(G)\rightarrow \{0,1,2\}$ defined by $f^{\prime}(v)=1$ and $f^{\prime}(b)=f(b)$ for all other $b\in V(G)$ would be a TRD-function on $G$ with $\omega(f^{\prime})=\omega(f)+1$, contradicting $uv$ being TRD-supercritical. Hence $g:V(G)\rightarrow \{0,1,2\}$ defined by $g(w)=g(v)=1$ and $g(b)=f(b)$ for all other $b\in V(G)$ is a $\gamma_{tR}(G)$-function. However, since $u$ is not isolated in $G[V_f^+]$, $g$ is also a TRD-function on $G-uw$, contradicting $uw$ being TRD-ER-critical. Hence $f(u)=0$ and $f(v)=2$. 

Now, $f(N_G(u))=0$, otherwise $f^{\prime }:V(G)\rightarrow \{0,1,2\}$ defined by $f^{\prime}(u)=1$ and $f^{\prime}(b)=f(b)$ for all other $b\in V(G)$ would be a TRD-function on $G$ with $\omega(f^{\prime})=\omega(f)+1$, contradicting $uv$ being TRD-supercritical. Furthermore, since $\text{deg}_G(u)\geq2$, there exists some vertex $y\in N_G(u)-\{w\}$. Note that $f(y)=0$. Hence $f^{\prime}:V(G)\rightarrow \{0,1,2\}$ defined by $f^{\prime}(y)=f^{\prime}(u)=1$ and $f^{\prime}(b)=f(b)$ for all other $b\in V(G)$ is a $\gamma_{tR}(G)$-function. However, $f^{\prime}$ is also a TRD-function on $G-uw$, contradicting $uw$ being TRD-ER-critical. Therefore $\text{deg}_G(u)=1$; that is, $uw$ is a pendant edge.~$\square $

\begin{cor}
\thlabel{super stable} If $G$ is a $\gamma_{tR}$-edge-supercritical graph with $\delta(G)\geq2$, then $G$ is $\gamma_{tR}$-ER-stable.
\end{cor}

\section{Future Work}
\label{Sec: Future}

Consider for a moment connected $6$-$\gamma_{tR}$-edge-supercritical graphs. We showed in Section \ref{Sec: 6-super} that, for any connected $6$-$\gamma_{tR}$-edge-supercritical graph $G$, $2\leq\text{diam}(G)\leq3$. Furthermore, note that each graph $G_r$, with $r\geq2$, introduced in Section \ref{Sec: 6-super} has diameter $3$. We now consider the $6$-$\gamma_{tR}$-edge-supercritical graphs $G$ for which $\text{diam}(G)=2$. We begin with the following lemma, which provides a lower bound for the minimum degree of a connected graph $G$ with diameter $2$, based on its TRD-number.

\begin{lemma}
\thlabel{diam2 deg bounds} If $G$ is a connected graph with $\emph{\text{diam}}(G)=2$ and $\gamma_{tR}(G)=k$, then $\delta(G)\geq\lfloor\frac{k}{2}\rfloor$.
\end{lemma}

\noindent \emph{Proof.} Suppose for a contradiction that there is a vertex $v\in V(G)$ with $\text{deg}(v)<\lfloor\frac{k}{2}\rfloor$. Since $\text{diam}(G)=2$, $N_G(v)$ is a dominating set of $G$. Thus the function $f:V(G)\rightarrow \{0,1,2\}$ defined by $f(v)=1$, $f(x)=2$ for all $x\in N_G(v)$, and $f(z)=0$ for all other $z\in V(G)$ is a TRD-function on $G$ with $\omega(f)\leq2(\lfloor\frac{k}{2}\rfloor-1)+1$. That is, $\omega(f)\leq2\lfloor\frac{k}{2}\rfloor-1<k$, contradicting $\gamma_{tR}(G)=k$.~$\square $

\begin{cor}
\thlabel{diam2 deg>2} If $G$ is a connected $\gamma_{tR}$-edge-supercritical graph with $\emph{\text{diam}}(G)=2$, then $\delta(G)\geq3$. 
\end{cor}

\vspace{0mm}

The previous corollary follows directly from \thref{no 5-super}. In light of this result, we present the following proposition which provides necessary conditions for a connected graph $G$ with $\text{diam}(G)=2$ to be $6$-$\gamma_{tR}$-edge-supercritical. Characterizing connected $6$-$\gamma_{tR}$-edge-supercritical graphs with diameter $2$, and indeed with diameter $3$, remain open problems.

\begin{lemma}
\thlabel{6 diam2 D=TD} If $G$ is a connected $6$-$\gamma_{tR}$-edge-supercritical graph with $\emph{\text{diam}}(G)=2$, then $G$ is $3$-$\gamma_{t}$-edge-critical and $3$-$\gamma$-edge-critical.
\end{lemma}

\noindent \emph{Proof.} Let $G$ be a connected $6$-$\gamma_{tR}$-edge-supercritical graph with $\text{diam}(G)=2$. Then, for any edge $e\in E(\overline{G})$, $\gamma_{tR}(G+e)=4$. Thus, by \thref{tR=34}, $\gamma_{t}(G+e)=\gamma(G+e)=2$. Now, \thref{tR=34} also implies that $\gamma_{t}(G)>2$. Furthermore, by \thref{t bounds}, $\gamma_{t}(G)\leq4$. If $\gamma_{t}(G)=4$, then $G$ is $4$-$\gamma_{t}$-edge-supercritical, which, since $G$ is connected, contradicts \thref{t super}. Hence $\gamma_{t}(G)=3$. Now, by \thref{TD vs D}, $2\leq\gamma(G)\leq3$. Suppose for a contradiction that $\gamma(G)=2$, and consider a $\gamma(G)$-set $S=\{u,v\}$. Note that, since $\gamma_{t}(G)=3$, $uv\in E(\overline{G})$. However, since $\text{diam}(G)=2$, there exists some $w\in N_G(u)\cap N_G(v)$. Hence $T=\{u,v,w\}$ is a $\gamma_{t}(G)$-set. But then $S\subset T$, contradicting \thref{tR=5}. Hence $\gamma(G)=3$, and thus $G$ is $3$-$\gamma_{t}$-edge-critical and $3$-$\gamma$-edge-critical.~$\square $ 

\begin{que}
Do there exist connected $6$-$\gamma_{tR}$-edge-supercritical graphs with diameter $2$?
\end{que}

\vspace{0mm}

Having demonstrated the existence of connected $6$-$\gamma_{tR}$-edge-supercritical graphs with diameter $3$ in Section \ref{Sec: 6-super}, we now consider the $\gamma_{tR}$-functions on these graphs $G_r$, where $r\geq2$.

\begin{prop}
\thlabel{G_r V+} For $r\geq2$, if $v\in V(G_r)$, then there exists a $\gamma_{tR}(G)$-function $f$ such that $v\in V_f^+$.
\end{prop}

\noindent \emph{Proof.} Let $v\in V(G_r)$, where $r\geq2$. Then, without loss of generality, $v\in\{x,y,z,u_0,u_1,w_0,w_1\}$. If $v\in\{x,y,z\}$ consider the function $f:V(G_r)\rightarrow\{0,1,2\}$ defined by $f(x)=f(y)=f(z)=2$ and $f(b)=0$ for all other $b\in V(G_r)$. Otherwise, if $v\in\{u_1,w_1\}$, consider the function $f:V(G_r)\rightarrow\{0,1,2\}$ defined by $f(u_1)=f(w_1)=f(z)=2$ and $f(b)=0$ for all other $b\in V(G_r)$. Otherwise, if $v=u_0$, consider the function $f:V(G_r)\rightarrow\{0,1,2\}$ defined by $f(u_0)=f(u_1)=f(w_2)=2$ and $f(b)=0$ for all other $b\in V(G_r)$. Otherwise, if $v=w_0$, consider the function $f:V(G_r)\rightarrow\{0,1,2\}$ defined by $f(x)=f(y)=f(w_0)=2$ and $f(b)=0$ for all other $b\in V(G_r)$. In any case, we have a $\gamma_{tR}(G)$-function $f$ such that $v\in V_f^+$, as required.~$\square $

\begin{cor}
\thlabel{G_r K_n crit} For $r\geq2$ and $n\geq3$, $G_r\cup K_n$ is $9$-$\gamma_{tR}$-edge-critical. 
\end{cor}

\noindent \emph{Proof.} Consider the graph $H\cong G_r\cup K_n$ where $r\geq2$ and $n\geq3$. Clearly $\gamma_{tR}(H)=9$. Consider an edge $uv\in E(\overline{H})$. If $uv\in E(\overline{G_r})$, \thref{G_r super} implies that $uv$ is supercritical, and thus critical, with respect to total Roman domination. Otherwise, suppose that $u\in V(G_r)$ and $v\in V(K_n)$. By \thref{G_r V+}, there exists a $\gamma_{tR}(G_r)$-function $g$ such that $u\in V_g^+$. Consider the function $f:V(G_r)\rightarrow\{0,1,2\}$ defined by $f(w)=g(w)$ for all $w\in V(G_r)$, $f(v)=2$, and $f(x)=0$ for all other $x\in V(K_n)$. Then $f$ is a TRD-function on $H+uv$ with $\omega(f)=8$, and hence $H$ is $9$-$\gamma_{tR}$-edge-critical.~$\square $

\vspace{0mm}

By \thref{no 5-super,disconnect 6-super}, the disjoint union of a disconnected $6$-$\gamma_{tR}$-edge-supercritical graph $G$ and $K_n$ for $n\geq3$ is $\gamma_{tR}$-edge-supercritical, and thus $\gamma_{tR}$-edge-critical. Moreover, it can be easily seen that the union of $\text{cor}(K_m)$ and $K_n$, with $m\geq4$ and $n\geq3$, is also $\gamma_{tR}$-edge-critical. In light of our previous result, we pose the following conjectures. Note that the second conjecture would be a direct result of the first.

\begin{conj}
\thlabel{join super} If $G$ is a $\gamma_{tR}$-edge-supercritical graph and $v\in V(G)$, then there exists a $\gamma_{tR}(G)$-function $f$ such that $v\in V_f^+$.  
\end{conj} 

\begin{conj}
\thlabel{join super} If $G$ is a $k$-$\gamma_{tR}$-edge-supercritical graph, then $G\cup K_n$ is $(k+3)$-$\gamma_{tR}$-edge-critical, for $n\geq3$.  
\end{conj}

\pagebreak
 
\label{Refs}

\end{document}